\documentclass[graybox]{svmult}

\usepackage{type1cm}        % activate if the above 3 fonts are
\usepackage{makeidx}         % allows index generation
\usepackage{graphicx}        % standard LaTeX graphics tool
\usepackage{multicol}        % used for the two-column index
\usepackage[bottom]{footmisc}% places footnotes at page bottom

\usepackage{mathtools}
\usepackage{newtxtext}       %
\usepackage{newtxmath}       % selects Times Roman as basic font

\usepackage{mathrsfs}
\usepackage{amssymb}
\DeclareMathAlphabet{\mathpzc}{OT1}{pzc}{m}{it}
\usepackage[capitalise]{cleveref}
\usepackage{subfig}
		\usepackage{tikz}
		\usepackage{pgfplots}
		\pgfplotsset{compat=newest}
		\usepackage{xcolor}
\definecolor{MPIblue}{RGB}{51,165,195}
\colorlet{mpiblue}{MPIblue}
\definecolor{MPIgrey}{RGB}{135,135,141}
\colorlet{mpigrey}{MPIgrey}
\definecolor{MPIgreen}{RGB}{0,118,117}
\colorlet{mpigreen}{MPIgreen}
\definecolor{MPIred}{RGB}{120,0,75}
\colorlet{mpired}{MPIred}
\definecolor{MPIsand}{RGB}{236,233,212}
\colorlet{mpisand}{MPIsand}
\definecolor{MPItext}{RGB}{56,60,60}
\colorlet{mpitext}{MPItext}
\definecolor{MPIbluetext}{RGB}{23,161,193}
\colorlet{mpibluetext}{MPIbluetext}
\definecolor{MPIbackground}{RGB}{239,238,227}
\colorlet{alg1col1}{MPIblue!30}
\colorlet{alg1col2}{MPIsand}

\usepackage{algorithm}
\usepackage{algorithmic}

\newcommand{\bsy}[1]{\ensuremath{\boldsymbol{{#1}}}}

\makeindex             % used for the subject index
\input{mymacros.sty}											

\usepackage{algorithm, algorithmic}
\begin{document}
\title*{Adaptive Interpolatory MOR by Learning the Error Estimator in the Parameter Domain}
\titlerunning{Adaptive Interpolatory MOR by Learning}
\author{Sridhar Chellappa, Lihong Feng, Valent\'in de la Rubia and Peter Benner}
\institute{Sridhar Chellappa \at Max Planck Institute for Dynamics of Complex Technical Systems, Sandtorstra\ss e 1, Magdeburg, 39106, Germany \email{chellappa@mpi-magdeburg.mpg.de}
\and Lihong Feng \at Max Planck Institute for Dynamics of Complex Technical Systems, Sandtorstra\ss e 1, Magdeburg, 39106, Germany \email{feng@mpi-magdeburg.mpg.de}
\and Valent\'in de la Rubia \at Departamento de Matem\'atica Aplicada a las TIC, ETSI de Telecomunicaci\'on, Universidad Polit\'ecnica de Madrid, 28040 Madrid, Spain \email{valentin.delarubia@upm.es}
\and Peter Benner \at Max Planck Institute for Dynamics of Complex Technical Systems, Sandtorstra\ss e 1, Magdeburg, 39106, Germany \email{benner@mpi-magdeburg.mpg.de}}
\maketitle

\abstract*{Interpolatory methods offer a powerful framework for generating reduced-order models (ROMs) for non-parametric or parametric systems with time-varying inputs. Choosing the interpolation points adaptively remains an area of active interest. A greedy framework has been introduced in \cite{morFenAB17,morFenB19} to choose interpolation points automatically using \textit{a posteriori} error estimators. Nevertheless, when the parameter range is large or if the parameter space dimension is larger than two, the greedy algorithm may take considerable time, since the training set needs to include a considerable number of parameters. As a remedy, we introduce an adaptive training technique by learning an efficient \textit{a posteriori} error estimator over the parameter domain. A fast learning process is created by interpolating the error estimator using radial basis functions (RBF) over a fine parameter training set, representing the whole parameter domain. The error estimator is evaluated only on a coarse training set including a few parameter samples. The algorithm is an extension of the work in~ \cite{morCheFB19} to interpolatory model order reduction (MOR) in frequency domain. Beyond the work in~\cite{morCheFB19}, we use a newly proposed inf-sup-constant-free error estimator in the frequency domain~\cite{morFenB19}, which is often much tighter than the error estimator using the inf-sup constant.  
%Two training sets are considered: one coarsely sampled and the other finely sampled. The error estimator is only evaluated over the former while a RBF surrogate model of the estimator is evaluated over the latter. The coarse training set is updated at each iteration with new parameter(s) from the fine training set which have the maximum interpolated error value. Such an approach overcomes any limitation that may arise from the large dimension or range of the parameter training set. 
Three numerical examples demonstrate the efficiency and validity of the proposed approach.}

\abstract{Interpolatory methods offer a powerful framework for generating reduced-order models (ROMs) for non-parametric or parametric systems with time-varying inputs. Choosing the interpolation points adaptively remains an area of active interest. A greedy framework has been introduced in \cite{morFenAB17,morFenB19} to choose interpolation points automatically using \textit{a posteriori} error estimators. Nevertheless, when the parameter range is large or if the parameter space dimension is larger than two, the greedy algorithm may take considerable time, since the training set needs to include a considerable number of parameters. As a remedy, we introduce an adaptive training technique by learning an efficient \textit{a posteriori} error estimator over the parameter domain. A fast learning process is created by interpolating the error estimator using radial basis functions (RBF) over a fine parameter training set, representing the whole parameter domain. The error estimator is evaluated only on a coarse training set including a few parameter samples. The algorithm is an extension of the work in~ \cite{morCheFB19} to interpolatory model order reduction (MOR) in frequency domain. Beyond the work in~\cite{morCheFB19}, we use a newly proposed inf-sup-constant-free error estimator in the frequency domain~\cite{morFenB19}, which is often much tighter than the error estimator using the inf-sup constant.  
%Two training sets are considered: one coarsely sampled and the other finely sampled. The error estimator is only evaluated over the former while a RBF surrogate model of the estimator is evaluated over the latter. The coarse training set is updated at each iteration with new parameter(s) from the fine training set which have the maximum interpolated error value. Such an approach overcomes any limitation that may arise from the large dimension or range of the parameter training set. 
Three numerical examples demonstrate the efficiency and validity of the proposed approach.}

\section{Introduction}
\label{sec:1}

MOR based on system theory and interpolation~\cite{morBauBBetal11, morFenAB17, morFenB14, morGri97, morGugAB08} has been developed as a class of efficient MOR 
methods among others. Detailed summary of those methods and comparison of them with other classes of MOR methods can be found in some survey papers and books~\cite{morAnt05, morAntSG01, morBauBF14, morBenGW15, morBenMS05, morGugSW13, morSchVR08}. A major advantage of the interpolatory methods is their flexibility in reducing systems with time or parameter varying inputs, since they are based on the transfer function or the input-output relation of the systems, which is independent of the input signal. On the contrary, the snapshot MOR methods, such as proper orthogonal decomposition (POD) and the reduced basis methods (RBM) are input-dependent, and are often less efficient in reducing systems with varying inputs as compared with the interpolatory MOR methods~\cite{morBenGW15}. 

A major topic of interest in interpolatory MOR methods is how to determine the interpolation points, so as to adaptively 
construct the ROM. Many methods have appeared in the last ten years, some are heuristic~\cite{morBecRK04, morFenKB15, morHetTF13, morLeeCF06}, some entail high computational complexity~\cite{morFenAB17, morWolPL11}, and some are inefficient for systems with more than one parameter~\cite{morGugAB08, morHunMS18}. Random interpolation points are used in~\cite{morBenGP19}.

Recently, a new error estimator for the reduced transfer function error and an algorithm for iteratively choosing the interpolation points are proposed in~\cite{morFenB19}, which overcomes many difficulties being faced by the above mentioned interpolatory methods. It is neither heuristic nor needs  a high computational cost. Moreover, it is a parametric MOR method and applicable to systems with more than two parameters. One shortcoming of the method is that the interpolation points are selected from a given training set, which must be decided a priori and becomes larger and larger with the increase of the parameter range or the parameter space dimension. Such a technique is standard also for the RBM, where a training set must be given before a greedy algorithm starts. This makes the greedy algorithm slow down when there is a large number of samples in the training set due to the large dimension or large range of the parameter domain. This is due to the fact that at each iteration of the greedy algorithm, an error estimator needs to be repeatedly computed for all the samples in the training set.  Many adaptive training techniques have been proposed recently for RBM~\cite{morCheFB19, morHaaDO11, morHesSZ14, morMadS13}. In contrast, no efficient training techniques are proposed for the interpolatory MOR methods, though similar greedy algorithms using fixed training sets are proposed in~\cite{morFenAB17, morFenB19}. In this work, we extend the adaptive training technique in~\cite{morCheFB19} for RBM to an adaptive training technique for the interpolatory MOR methods in~\cite{morFenAB17, morFenB19}. 

The main contribution of this work is an efficient algorithm to adaptively choose interpolation points for parametric, linear time-invariant (LTI) systems having a wide range of parameter values or with a large parameter space dimension. Compared with the greedy algorithms proposed in~\cite{morFenAB17, morFenB19}, we have added two new ingredients to the greedy algorithms: (i) a sharp \textit{a posteriori} output (or state) error estimator with low computational costs and (ii) a surrogate with even lower computational costs for learning the error estimator over the whole parameter domain. 
The aim is instead of computing the error estimator over the whole parameter domain the for ROM construction, a surrogate estimator is computed. In this way, the error estimator is computed only on a coarse training set at each iteration of the greedy algorithm, and for parameters outside of the coarse training set, the surrogate estimator is computed. Finally,  the training set needs to be initialized by including only a few parameter samples, and can be iteratively updated using the surrogate estimator instead of the error estimator itself. 
As a consequence, a significant amount of computational cost can be saved for such systems. 

The idea is similar to the one in~\cite{morCheFB19} for the RBM. However, in~\cite{morCheFB19}, an error estimator in time domain is used, where 
the inf-sup constant needs to be computed for each parameter in the training set, which is computationally inefficient for large-scale systems. In this work, we use an inf-sup-constant-free error estimator newly proposed in~\cite{morFenB19}. It is suitable for interpolatory MOR methods, since it estimates the transfer function error in the frequency domain.  

The paper is organized as follows. In \cref{sec:2} we briefly review interpolatory MOR methods based on projection. The greedy interpolatory methods~\cite{morFenAB17, morFenB19} for parametric systems are reviewed in~\cref{sec:greedy}. In~\cref{sec:3}, we introduce the basic idea of RBF interpolation and elaborate on the process of learning the error estimator using a surrogate estimator constructed by RBF interpolation. Based on this surrogate estimator, we propose the greedy algorithm \textsf{IPSUE} with adaptive training technique for adaptively choosing the interpolation points in a more efficient and fully adaptive way. We present numerical results on three real-word examples in Section 4 to show the robustness of \textsf{IPSUE} and conclude the work in the end.

\section{Interpolatory MOR}
\label{sec:2}
In this work, we are interested in MOR of parametric LTI systems in the \textit{state-space representation} given by
		\begin{equation}
				\Sigma(\bsy{\mu}):\begin{cases}
							\begin{aligned}
							{\bf E}\dot{\bx}(t,\bsy{\mu}) &= \bA(\bsy{\mu}) \bx(t,\bsy{\mu}) + \bB(\bsy{\mu}) \bu(t),\\
							\by(t,\bsy{\mu}) &= \bC(\bsy{\mu}) \bx(t,\bsy{\mu}),\,\, \bx(0,\bsy{\mu}) = 0.
							\end{aligned}
						\end{cases}
			\label{eq:LTI}
		\end{equation}
Here, $\bsy{\mu} := \big[ \mu^{1}, \mu^{2}, \ldots, \mu^{d}\big]^{\textsf{T}} \in \R^{d}$ is the vector of parameters (geometric or physical). $\bx(t,\bsy{\mu}) \in \R^{n}$ is the state vector and $n$ is typically very large. $\bu(t) \in \R^{m}$ is the input vector and $\by(t,\bsy{\mu}) \in \R^{p}$ is the output vector. $\bA(\bsy{\mu}) \in \R^{n \times n}$ is the state matrix, $\bB(\bsy{\mu}) \in \R^{n \times m}$ is the input matrix, $\bC(\bsy{\mu}) \in \R^{p \times n}$ is the output matrix. For the case when $m = p =1$, \cref{eq:LTI} is referred to as a single-input, single-output (SISO) system. Otherwise, it is known as multi-input, multi-output (MIMO) system. 

The ROM we seek should preserve the same structure of the FOM but have a much smaller dimension. We assume that the state vector lies (approximately) in the span of a low-dimensional linear subspace $\mathcal{V} \subset \R^{n \times r}$, $r \ll n$, such that, $\bx(t,\bsy{\mu}) \approx \bV \hat{\bx}(t,\bsy{\mu})$. Column vectors in the matrix $\bV \in \R^{n \times r}$ constitute an orthogonal basis of $\mathcal{V}$. Replacing $\bx(t,\bsy{\mu})$ in~\cref{eq:LTI} with its approximation $\bV \hat{\bx}(t,\bsy{\mu})$ and further imposing Petrov-Galerkin projection on the residual introduced by the approximation in a test subspace $\mathcal{W} \subset \R^{n \times r}$ leads to
	\begin{equation*}
		\bW^{\textsf{T}} \bigg( \bV \dot{\hat{\bx}}(t,\bsy{\mu}) -  \bA(\bsy{\mu}) \bV \hat{\bx}(t,\bsy{\mu}) - \bB(\bsy{\mu}) \bu(t)\bigg) \equiv \mathbf{0},
	\end{equation*}
where, the column vectors in the matrix $\bW \in \R^{n \times r}$ correspond to an orthogonal basis of $\mathcal{W}$. The resulting ROM is given as
		\begin{equation}
			\hat{\Sigma}(\bsy{\mu}):\begin{cases}
			\begin{aligned}
			\hat {\bf E}\dot{\hat{\bx}}(t,\bsy{\mu}) &= \hat{\bA}(\bsy{\mu}) \hat{\bx}(t,\bsy{\mu}) + \hat{\bB}(\bsy{\mu}) \bu(t),\\
			\hat{\by}(t,\bsy{\mu}) &= \hat{\bC}(\bsy{\mu}) \hat{\bx}(t,\bsy{\mu}),\,\, \hat{\bx}(0,\bsy{\mu}) = 0.
			\end{aligned}
			\end{cases}
		\label{eq:LTI_rom}
		\end{equation}

Here, $\hat{\bx}(t,\bsy{\mu}) \in \R^{r}$ is the reduced state vector, $\hat{\bf E}(\bsy{\mu})=\bW^{\textsf{T}} {\bf E} (\bsy{\mu}) \bV \in \R^{r \times r}$, $\hat{\bA}(\bsy{\mu})=\bW^{\textsf{T}} \bA(\bsy{\mu}) \bV \in \R^{r \times r}, \hat{\bB}(\bsy{\mu})=\bW^{\textsf{T}} \bB(\bsy{\mu}) \in \R^{r \times m}, \hat{\bC}(\bsy{\mu})= \bC(\bsy{\mu}) \bV \in \R^{p \times r}$ are the reduced system matrices, and $\hat{\by}(t,\bsy{\mu})$ is the reduced output vector. 
The goal of MOR is to find the two subspaces $\mathcal{V}, \mathcal{W} \in \R^{n \times r}$. Different MOR methods vary in how they generate the matrices $\bW, \bV$.

Interpolatory MOR methods construct $\bV, \bW \in \R^{n \times r}$ based on the transfer function of the system, which is independent of the input signal. The transfer function of the system described in \cref{eq:LTI} is given by
	\begin{equation}
		\bH(\tilde{\bsy{\mu}}) :=  \bC(\bsy{\mu}) \big( \overbrace{s {\bf E} - \bA(\bsy{\mu})}^{=: \mathscr{A}(\tilde{\bsy{\mu}})} \big)^{-1} \bB(\bsy{\mu}).
		\label{eq:LTI_s}
	\end{equation}
Here, $\tilde{\bsy{\mu}} := \big[s, \mu^{1}, \mu^{2}, \ldots, \mu^{d}\big]^{\textsf{T}} \in \R^{d+1}$ is the vector of parameters with the additional Laplace variable $s \in \jmath \R^{}$, where $\jmath$ is the imaginary unit. The corresponding ROM of \cref{eq:LTI_s} is of the form,
	\begin{equation}
		\hat{\bH}(\tilde{\bsy{\mu}})  := \hat{\bC}(\tilde{\bsy{\mu}}) \hat{\mathscr{A}}(\tilde{\bsy{\mu}})^{-1} \hat{\bB}(\bsy{\mu}),
	\end{equation}
with $\hat{\mathscr{A}}(\tilde{\bsy{\mu}}) := s \hat {\bf E} - \hat{\bA}(\bsy{\mu})$.

Many interpolatory methods have been proposed for linear systems, especially for linear non-parametric systems. The most representative methods are the moment-matching methods~\cite{morGri97, morGugAB08},
where the $\mathcal H_2$-optimal method IRKA~\cite{morGugAB08} constructs a ROM satisfying the necessary conditions of local optimality. All these methods are known to be applicable to non-parametric systems. Later, IRKA is extended to MOR for  parametric systems~\cite{morBauBBetal11}, where some pairs of projection matrices are constructed for given samples of parameters, then they are combined together to get the final pair of projection matrices. No rule is used for selecting the samples. In~\cite{morHunMS18}, a method for parametric systems is proposed based on $\mathcal H_2\otimes \mathcal L_2$-optimality, but is only applicable to systems with one parameter and is facing high computational complexity for systems with $n \geq 1000$. 

Choosing interpolation points using a greedy algorithm guided by an \textit{a posteriori} error bound is proposed in~\cite{morFenAB17}. However, computing the error estimator needs to compute the smallest singular values of the large matrix $\mathscr{A}(\tilde{\bsy{\mu}})$, the inf-sup constant.
An inf-sup-constant-free error estimator is newly proposed in~\cite{morFenB19}, which can be efficiently computed, and is also much tighter than the error bound~\cite{morFenAB17} for many systems with small inf-sup constants. A similar greedy algorithm is proposed in~\cite{morFenB19} for choosing the interpolation points using the new error estimator. The adaptive training approach proposed in this work is based on the greedy algorithm and the new error estimator in~\cite{morFenB19}. In the next section, we briefly review the error estimator and the corresponding greedy algorithm.

\section{Greedy Method for Choosing Interpolation Points}
\label{sec:greedy}
The transfer function can be seen as a mapping from the space of inputs $\R^{m}$ to the space of outputs $\R^{p}$ passing through a high dimensional intermediate state in $\R^{n}$. 
If we look at the matrix product $\mathscr{A}^{-1}(\tilde{\bsy{\mu}})\bB(\bsy{\mu})$ in ${\bf H}(\bsy{\mu})$,  we may consider the primal system,
	\begin{equation}
		\begin{aligned}
			\mathscr{A}(\tilde{\bsy{\mu}}) \bX_{\text{pr}}(\tilde{\bsy{\mu}}) &= \bB(\bsy{\mu}). \label{eq:primal_state}
		\end{aligned}
	\end{equation}
Here, $\bX_{\text{pr}}(\tilde{\bsy{\mu}}) \in \R^{n}$ is the primal state vector. The reduced primal system is defined as,
	\begin{equation}
		\begin{aligned}
			\hat{\mathscr{A}}(\tilde{\bsy{\mu}}) \hat{\bX}_{\text{pr}}(\tilde{\bsy{\mu}}) &= \hat{\bB}(\bsy{\mu}). \label{eq:primal_state_rom}
		\end{aligned}
	\end{equation}
The approximate primal solution is given by $\tilde{\bX}_{\text{pr}} (\tilde{\bsy{\mu}}) := \bV \hat{\bX}_{\text{pr}}(\tilde{\bsy{\mu}})$ and the corresponding residual is
	\begin{equation}
		\br_{\text{pr}}(\bsy{\mu}) = \bB(\bsy{\mu}) - \mathscr{A}(\tilde{\bsy{\mu}}) \tilde{\bX}_{\text{pr}}(\tilde{\bsy{\mu}}).
		\label{eq:primal_residual}
	\end{equation}
Additionally, by considering the matrix product $\bC(\bsy{\mu})\mathscr{A}^{-1}(\tilde{\bsy{\mu}})$ in ${\bf H}(\bsy{\mu})$, we have the following dual system,
	\begin{equation}
		\begin{aligned}
			\mathscr{A}^{\textsf{T}}(\tilde{\bsy{\mu}}) \bX_{\text{du}}(\tilde{\bsy{\mu}}) &= \bC^{\textsf{T}}(\bsy{\mu}). \label{eq:dual_state}
		\end{aligned}
	\end{equation}
$\bX_{\text{du}}(\tilde{\bsy{\mu}}) \in \R^{n}$ is the dual state vector. The reduced dual system is given as,
	\begin{equation}
		\begin{aligned}
			\hat{\mathscr{A}}^{\textsf{T}}(\tilde{\bsy{\mu}}) \hat{\bX}_{\text{du}}(\tilde{\bsy{\mu}}) &= \hat{\bC}^{\textsf{T}}(\bsy{\mu}). \label{eq:dual_state_rom}
		\end{aligned}
	\end{equation}
The approximate dual solution is given by $\tilde{\bX}_{\text{du}} (\tilde{\bsy{\mu}}) := \bV_{\text{du}} \hat{\bX}_{\text{du}}(\tilde{\bsy{\mu}})$ and the corresponding residual is,
	\begin{equation}
		\br_{\text{du}}(\bsy{\mu}) = \bC^{\textsf{T}}(\bsy{\mu}) - \mathscr{A}^{\textsf{T}}(\tilde{\bsy{\mu}}) \tilde{\bX}_{\text{du}}(\tilde{\bsy{\mu}}).
		\label{eq:dual_residual}
	\end{equation}	

For parametric LTI systems, adopting the spirit of the RBM, \cite{morFenAB17} introduced a method to automatically generate a ROM through a greedy algorithm. The authors introduce a primal-dual residual-based \textit{a posteriori} error estimator for the transfer function approximation error $\|\bH(\tilde{\bsy{\mu}}) - \hat{\bH}(\tilde{\bsy{\mu}})\|$, for both SISO and MIMO systems. For SISO systems it reads,
		\begin{equation}
			| \bH(\tilde{\bsy{\mu}}) - \hat{\bH}(\tilde{\bsy{\mu}}) | \leq \frac{\| \br_{\text{pr}}(\tilde{\bsy{\mu}}) \|_{2}  \|\br_{\text{du}}(\tilde{\bsy{\mu}})\|_{2} }{\sigma_{\text{min}}(\bsy{\mu})}.
			\label{eq:error_estm}
		\end{equation}
Here, $\sigma_{\text{min}}(\bsy{\mu})$, called the inf-sup constant, is the smallest singular value of the matrix $\mathscr{A}(\bsy{\mu})$ as defined in \cref{eq:LTI_s}. The primal and dual residuals are given in \cref{eq:primal_residual} and \cref{eq:dual_residual}, respectively. The work \cite{morFenB19} improves the method in \cite{morFenAB17} by avoiding the calculation of the inf-sup constant required for the error estimator. This is achieved by introducing a dual-residual system,
		\begin{equation}
			\mathscr{A}^{\textsf{T}}(\tilde{\bsy{\mu}}) \mathbf{e}_{\text{du}}(\tilde{\bsy{\mu}}) = \br_{\text{du}}(\tilde{\bsy{\mu}}).
			\label{eq:aux_du_resd}
		\end{equation}
The following theorem from \cite{morFenB19} gives the \textit{a posteriori} error bound,
	\begin{proposition}
		The transfer function approximation error can be bounded as,
			\begin{equation*}
			|\bH(\tilde{\bsy{\mu}}) - \hat{\bH}(\tilde{\bsy{\mu}})| \leq \big| \tilde{\bX}_{\text{du}}^{\textsf{T}}(\tilde{\bsy{\mu}}) \br_{\text{pr}}(\tilde{\bsy{\mu}}) \big| + \big|\mathbf{e}_{\text{du}}^{\textsf{T}}(\tilde{\bsy{\mu}}) \br_{\text{pr}}(\tilde{\bsy{\mu}})\big|.
			\end{equation*}
		\label{thm:2}
	\end{proposition}
For a proof of \cref{thm:2}, we refer to \cite{morFenB19}. In this form, the error bound is not computationally efficient since one needs to solve the full order dual-residual system~(\ref{eq:aux_du_resd}) to obtain $\mathbf{e}_{\text{du}}(\tilde{\bsy{\mu}})$. 
Instead, system~(\ref{eq:aux_du_resd}) is reduced by an orthogonal matrix $\bV_{e} \in \R^{n \times \ell}$ as below,
		\begin{equation}
			\hat{\mathscr{A}}_{e}^{\textsf{T}}(\tilde{\bsy{\mu}}) \hat{\mathbf{e}}_{\text{du}}(\tilde{\bsy{\mu}}) = \hat{\br}_{\text{du,e}}(\tilde{\bsy{\mu}}),
			\label{eq:aux_du_resd_rom}
		\end{equation}
where $\hat{\mathscr{A}}_{e} := \bV_{e}^{\textsf{T}}\, \mathscr{A}(\tilde{\bsy{\mu}})\, \bV_{e}$ and $\hat{\br}_{\text{du,e}} := \bV_{e}^{\textsf{T}}\, \br_{\text{du}}(\tilde{\bsy{\mu}})$. 
The projection matrices $\bV, \bV_{\text{du}}$ and $\bV_{e}$ corresponding to the primal, dual and the dual-residual system are generated offline. 

By using the approximate solution to the dual-residual system, an efficiently computable error estimator is obtained.
	\begin{equation}
			|\bH(\tilde{\bsy{\mu}}) - \hat{\bH}(\tilde{\bsy{\mu}})| \lesssim \big| \tilde{\bX}_{\text{du}}^{\textsf{T}}(\tilde{\bsy{\mu}}) \br_{\text{pr}}(\tilde{\bsy{\mu}}) \big| + \big|\tilde{\mathbf{e}}_{\text{du}}^{\textsf{T}}(\tilde{\bsy{\mu}}) \br_{\text{pr}}(\tilde{\bsy{\mu}})\big| =: \Delta(\tilde{\bsy{\mu}}),	
			\label{eq:err_est_computable}	
	\end{equation}
where $\tilde{\mathbf{e}}_{\text{du}}(\tilde{\bsy{\mu}}) := \bV_{e}\, \hat{\mathbf{e}}_{\text{du}}(\tilde{\bsy{\mu}})$.
For ease of comparison, we first present the greedy algorithm for parametric systems introduced in \cite{morFenB19} as~\cref{alg:1}. 
% Alg Feng/Benner
\renewcommand{\algorithmicrequire}{\textbf{Input:}}
\renewcommand{\algorithmicensure}{\textbf{Output:}}
\begin{algorithm}[t!]
	\caption{\textsf{Greedy ROM Construction for Parametric Systems} \cite{morFenB19}}
	\begin{algorithmic}[1]
		\label{alg:1}
		\REQUIRE System matrices $\bA(\bsy{\mu}), \bB(\bsy{\mu}), \bC(\bsy{\mu})\,$, Training set $\Xi$ of cardinality $N_{\mu}\,$ covering the interesting parameter ranges, Tolerance $\epsilon_{tol}$.
		
		\ENSURE Projection matrix $\bV$.
		
		\STATE Initialize $\bV = [\,]$, $\bV_{\text{du}} = [\,]$,  $\bV_{e} = [\,]$, $\epsilon = 1 + \epsilon_{tol}$, fix $\eta$, the number of moments to be matched.
		
		\STATE Initial interpolation point $\tilde{\bsy{\mu}}^{1}$: the first sample in $\Xi$. $\tilde{\bsy{\mu}}_{\alpha}^{1}$: the last sample in $\Xi$. Set $i = 1$.
		
		\WHILE{$\epsilon$ $>$ $\epsilon_{tol}$} 
		
		\STATE Solve \cref{eq:primal_state} at interpolation point $\tilde{\bsy{\mu}} = \tilde{\bsy{\mu}}^{i}$ and update projection matrix
		\begin{equation*}
		\bV = \textsf{orth} \big( [\bV\,\,\texttt{mmm}(\mathscr{A}(\tilde{\bsy{\mu}}^{i}), \bB(\tilde{\bsy{\mu}}^{i}), \eta, \tilde{\bsy{\mu}}^{i} )]\big).
		\end{equation*}
		\STATE Solve \cref{eq:dual_state} at interpolation point $\tilde{\bsy{\mu}} = \tilde{\bsy{\mu}}^{i}$ and update projection matrix
		\begin{equation*}
		\bV_{\text{du}} = \textsf{orth}\big([\bV_{\text{du}}\,\,\texttt{mmm}(\mathscr{A}^{\textsf{T}}(\tilde{\bsy{\mu}}^{i}), \bC^{\textsf{T}}(\tilde{\bsy{\mu}}^{i}), \eta, \tilde{\bsy{\mu}}^{i} )]\big).
		\end{equation*}
		\STATE Solve \cref{eq:aux_du_resd} at interpolation point $\tilde{\bsy{\mu}} = \tilde{\bsy{\mu}}_{\alpha}^{i}$ and update projection matrix,
		\begin{equation*}
			\bV_{e} = \textsf{orth}\big([\bV_{e}\,\,\bV_{\text{du}}\,\,\texttt{mmm}(\mathscr{A}^{\textsf{T}}(\tilde{\bsy{\mu}}_{\alpha}^{i}), \bC^{\textsf{T}}(\tilde{\bsy{\mu}}_{\alpha}^{i}), \eta, \tilde{\bsy{\mu}}_{\alpha}^{i} )]\big).
		\end{equation*}		
				
		\STATE $i = i + 1$. 
		
		\STATE $\tilde{\bsy{\mu}}^{i} = \arg \max \limits_{\tilde{\bsy{\mu}} \in \Xi}  \Delta(\tilde{\bsy{\mu}})$.
		
		\STATE $\tilde{\bsy{\mu}}_{\alpha}^{i} = \arg \max \limits_{\tilde{\bsy{\mu}} \in \Xi} \big|\tilde{\mathbf{e}}_{\text{du}}^{\textsf{T}}(\tilde{\bsy{\mu}}) \br_{\text{pr}}(\tilde{\bsy{\mu}})\big|$.		
		\STATE $\epsilon = \Delta(\tilde{\bsy{\mu}}^{i})$.	
		\ENDWHILE
	\end{algorithmic}
\end{algorithm}

It is automatic apart from the need for determining, \textit{a priori}, a suitable training set $\Xi$. The method proceeds by picking points from $\Xi$ that maximize the error estimator at every iteration and updating the three projection matrices $\bV, \bV_{\text{du}}, \bV_{e}$. However, there is no principled way to select the training set \textit{a priori}. If not adequately sampled, the training set may result in a ROM whose error is not uniformly below the tolerance. When the parameters involved can take on a wide range of values, or if many parameters are involved, then the number of parameter samples in $\Xi$ becomes large and the offline computation costs rise. We propose to solve this issue by constructing a surrogate model for $\Delta(\tilde{\bsy{\mu}})$ in \cref{eq:err_est_computable} and assure that computing the surrogate is much cheaper than computing the error estimator itself. The next section discusses this idea.	
\begin{remark}
	The algorithm is also applicable to MIMO systems. In this case the transfer function is matrix-valued. The key is how to compute the error estimator $\Delta(\tilde{\bsy \mu})$. We first estimate the error of the reduced transfer function entry-wise, i.e.
		\begin{equation}
			|\bH_{ij}(\tilde{\bsy{\mu}}) - \hat{\bH}_{ij}(\tilde{\bsy{\mu}})| \lesssim \big| \tilde{\bX}_{\text{du}}^{\textsf{T}}(\tilde{\bsy{\mu}}) \br_{\text{pr}}(\tilde{\bsy{\mu}}) \big| + \big|\tilde{\mathbf{e}}_{\text{du}}^{\textsf{T}}(\tilde{\bsy{\mu}}) \br_{\text{pr}}(\tilde{\bsy{\mu}})\big| =: \Delta_{ij}(\tilde{\bsy{\mu}}).		
			\label{eq:error_estm_MIMO}	
		\end{equation}
	Note that the $ij$-th entry of the transfer function corresponds to the input signal at the $j$-th input port and the signal at the $i$-th output port. 		
	Then, $\tilde{\bX}_{\text{du}}(\tilde{\bsy{\mu}})$ is the solution to the dual system by considering the right hand side as the $i$-th row vector of $\bC(\bsy{\mu})$, namely,  $\bC^{\textsf{T}}(:\,, i)$ in \cref{eq:dual_state_rom}. Correspondingly, the residual $\br_{\text{pr}}(\tilde{\bsy{\mu}})$ is obtained by solving \cref{eq:primal_state_rom} with the right hand side being $\bB(:\,,j)$, the $j$-th column of $\bB(\bsy{\mu})$. Then, $\Delta(\tilde{\bsy{\mu}}) = \arg \max \limits_{i,j} \Delta_{ij}(\tilde{\bsy{\mu}})$.
\end{remark}
\begin{remark}
	In order to build the projection matrices ($\bV, \bV_{\text{du}}, \bV_{e}$), \cite{morFenB19} makes use of the multi-moment matching (MMM) algorithm from \cite{morFenB14}. The algorithm provides an orthogonal basis for the solution at a given interpolation point, obtained through multivariate power series expansion of the state vector. To be focused on our main contribution, we refer to \cite{morFenB14,morFenB19} for detailed computations. For use in the proposed algorithm, we give below the call to the algorithm in MATLAB\textsuperscript \textregistered notation,
		\begin{equation*}
		\bV_{\textsf{mmm}} = \texttt{mmm}(\mathpzc{A}(\tilde{\mu}_{0}), \mathpzc
		{B}(\tilde{\mu}_{0}), \eta, \tilde{\mu}_{0}).		
		\end{equation*}
Here, $\mathpzc{A}(\tilde{\mu}_{0})$ denotes an arbitrary matrix evaluated at a given interpolation point $\tilde{\mu}_{0}$, $\mathpzc
{B}(\tilde{\mu}_{0})$ corresponds to the right hand side matrix in \cref{eq:primal_state,eq:dual_state,eq:aux_du_resd}, respectively. $\eta$ is the number of moments to be matched in the power series.	When $\eta = 0$, the MMM algorithm is equivalent to RBM, see~\cite{morFenB19} for more explanations. \end{remark}
\section{Adaptive Training by Learning the Error Estimator in the Parameter Domain}
\label{sec:3}

In this section, we propose an adaptive training technique, so that the greedy algorithm 
starts with a training set with small cardinality, which is then updated iteratively by using a surrogate error estimator. 
Different works have considered surrogate models of error estimators/indicators \cite{morCheFB19, morDroC15, morTaineA15}. All of these consider a surrogate in the context of the RBM. In this work, we deal with the frequency domain interpolatory MOR methods and focus on a surrogate model of an error estimator for the transfer function approximation error. The method we propose here is essentially an extension of the RBF-based error surrogate in~\cite{morCheFB19} to the frequency domain. Beyond the work in~\cite{morCheFB19}, we introduce a learning process in~\cref{subsec:learn} to show in detail how a surrogate estimator is constructed for any parameter in the whole parameter domain. We begin by introducing the method of RBF interpolation.
\subsection{Radial Basis Functions}
\label{subsec:rbf}
Radial Basis Functions belong to the family of \textit{kernel methods} and are a popular technique to generate surrogate models of multivariate functions $f : \R^{d} \mapsto \R^{}$, defined in a domain $\Omega \subset \R^{d}$. It may be the case that the function $f$ itself is unknown and one only knows a set of inputs $M = \{ \tilde{\bsy{\mu}}_{1}, \tilde{\bsy{\mu}}_{2}, \ldots, \tilde{\bsy{\mu}}_{\ell}\} \in \Omega$ and the corresponding function evaluations $F = \{f_{1}, f_{2}, \ldots, f_{\ell}\} \subset \R^{}$. Or, it may be the case that $f$ is known, but very expensive to evaluate repeatedly. For either case, RBF serves to generate an interpolant $g : \R^{d} \mapsto \R^{}$ of $f$ given by
		\begin{equation}
			g(\tilde{\bsy{\mu}}) = \sum_{i = 1}^{\ell} c_{i} \Phi(\| \tilde{\bsy{\mu}} - \tilde{\bsy{\mu}}_{i} \|), \, \forall \tilde{\bsy{\mu}} \in \Omega,
			\label{eq:rbf-interpolant}
		\end{equation}
such that it interpolates the original function at the set of input points (or centers) in $M$, i.e., $f(\tilde{\bsy{\mu}}_{i}) = g(\tilde{\bsy{\mu}}_{i}),\, i = 1, \ldots, \ell$. Moreover, $| f(\tilde{\bsy{\mu}}) - g(\tilde{\bsy{\mu}}) | \ll \texttt{tol}, \, \forall \tilde{\bsy{\mu}} \in \Omega$.
The functions $\Phi(\cdot)$ are the kernels defined as $\Phi(\tilde{\bsy{\mu}}_{1}, \tilde{\bsy{\mu}}_{2}) := \Phi(\| \tilde{\bsy{\mu}}_{1} - \tilde{\bsy{\mu}}_{2} \|), \forall \tilde{\bsy{\mu}}_{1}, \tilde{\bsy{\mu}}_{2} \in \Omega$. They are called radial basis functions owing to their radial dependence on $\tilde{\bsy{\mu}}$. The coefficients $\{c_{i}\}_{i = 1 }^{\ell}$ are determined by solving the linear system of equations,
		\begin{equation}
			\underbrace{
						\begin{bmatrix}
							\Phi(\tilde{\bsy{\mu}}_{1}, \tilde{\bsy{\mu}}_{1}) & \Phi(\tilde{\bsy{\mu}}_{1}, \tilde{\bsy{\mu}}_{2}) & \cdots & \Phi(\tilde{\bsy{\mu}}_{1}, \tilde{\bsy{\mu}}_{\ell}) \\
							\Phi(\tilde{\bsy{\mu}}_{2}, \tilde{\bsy{\mu}}_{1}) & \Phi(\tilde{\bsy{\mu}}_{2}, \tilde{\bsy{\mu}}_{2}) & \cdots & \Phi(\tilde{\bsy{\mu}}_{2}, \tilde{\bsy{\mu}}_{\ell}) \\
							\vdots & \vdots & \ddots & \vdots \\
							\Phi(\tilde{\bsy{\mu}}_{\ell}, \tilde{\bsy{\mu}}_{1}) & \Phi(\tilde{\bsy{\mu}}_{\ell}, \tilde{\bsy{\mu}}_{2}) & \cdots & \Phi(\tilde{\bsy{\mu}}_{\ell}, \tilde{\bsy{\mu}}_{\ell})
						\end{bmatrix}
					   }_{\mathbf{R}}
			\underbrace{
						\begin{bmatrix}
							c_{1} \\ c_{2} \\ \vdots \\ c_{\ell}
						\end{bmatrix}
					   }_{c}
					   =
			\underbrace{
						\begin{bmatrix}
							f(\tilde{\bsy{\mu}}_{1})\\ f(\tilde{\bsy{\mu}}_{2}) \\ \vdots \\ f(\tilde{\bsy{\mu}}_{\ell})
						\end{bmatrix}.
					   }_{\mathfrak{d}}
			\label{eq:rbfsystem}			
		\end{equation}
We need $\bR$ to be invertible. Assuming that the centers $\tilde{\bsy \mu}_{i}$ are pairwise distinct, it can be shown that $\mathbf{R}$ is positive definite for a suitable choice of the RBF $\Phi(\cdot)$ and thus \cref{eq:rbfsystem} has a unique solution. The class of RBF giving rise to positive definite $\bR$ is limited. As a workaround, some additional constraints are imposed in practice, i.e.
		\begin{equation*}
			\sum_{j=1}^{D} c_{i} p_{j}(\tilde{\bsy{\mu}}) = 0, \quad i = 1, 2, \ldots , \ell,
		\end{equation*}
so that a larger class of $\Phi(\cdot)$ can be admitted.		
The functions $p_{1}, p_{2}, \ldots, p_{D}$  are a basis of the polynomial space with suitable degree. In practice, we choose $D$ to be equal to the number of scalar parameters $d$. With the new conditions imposed, the radial basis interpolant now becomes,
		\begin{equation}
			g(\tilde{\bsy{\mu}}) := \sum_{i=1}^{\ell} c_{i} \Phi(\| \tilde{\bsy{\mu}} - \tilde{\bsy{\mu}}_{i} \|) + \sum_{j=1}^{D} \lambda_{j} p_{j}(\tilde{\bsy{\mu}}).
			\label{eq:rbf-interpolant-expanded}
		\end{equation}
This results in a saddle-point system of dimension $N_{\text{RBF}} := (D + \ell) \times (D + \ell)$,
		\begin{equation}
			\begin{bmatrix}
				\mathbf{R} & \mathbf{P}\\
				\mathbf{P}^{\textsf{T}} & 0
				\end{bmatrix}
				\begin{bmatrix}
				c \\ \lambda
				\end{bmatrix}
				=
				\begin{bmatrix}
				\mathfrak{d} \\ 0
			\end{bmatrix}.
		\label{eq:rbfsystem-saddle}
		\end{equation}
With a proper choice of $p_{1}, p_{2}, \ldots, p_{D}$, the augmented coefficient matrix is positive definite for a wider choice of kernel functions $\Phi(\cdot)$. 		
We refer to \cite{Wed05} for an exhaustive theoretical analysis of RBFs and the recent review paper \cite{SanH19} that analyses RBFs in the larger context of kernel based surrogate models.

\subsection{Learning the Error Estimator over the Parameter Domain}
\label{subsec:learn}
As highlighted in the Introduction, one of the main bottlenecks of the standard greedy algorithm is that the error estimator ($\Delta(\tilde{\bsy{\mu}})$) needs to be determined at every parameter in the training set. 

In order to evaluate it cheaply, we first construct a surrogate model of the error estimator by learning the error estimator in the whole parameter domain using RBF interpolation. We have the multivariate function $f(\cdot) := \Delta(\tilde{\bsy{\mu}})$ and the learning step involves determining the coefficients $c$ in \cref{eq:rbfsystem-saddle}.
First, we evaluate the error estimator at a small number of parameters in a coarse training set ($\Xi_{c} : [\tilde{\bsy{\mu}}_{1}, \ldots, \tilde{\bsy{\mu}}_{N_{c}}]$). These points shall serve as the centers $\tilde{\bsy \mu}$ of the RBF interpolation with $N_{c}$ as the number of centers. Note that, with regards to the discussion in \cref{subsec:rbf}, we have $\ell = N_{c}$.

Next, we define the kernel function ($\Phi(\cdot)$) and setup the linear system defined in \cref{eq:rbfsystem-saddle}. Many choices of the kernel function exist, and in the numerical experiments we have used the inverse multiquadric and the thin-plate spline kernel functions. For a deeper discussion, we refer to \cite{Wed05}. We note here that the assembling of the kernel matrix $\mathbf{R}$ can be done efficiently and software implementations exist to achieve this \cite{SanH19}. The right hand side is defined by $\mathfrak{d} := [\Delta(\tilde{\bsy{\mu}}_{1}), \ldots, \Delta(\tilde{\bsy{\mu}}_{N_{c}})]$. 

\cref{eq:rbfsystem-saddle} constitutes a small, dense system of linear equations. The computational cost of its evaluation scales as $\mathcal{O}((N_{c}+D)^{3})$. However, since $N_{c}, D$ are small, the cost remains under control.
Once knowing $c$ after solving \cref{eq:rbfsystem-saddle}, the interpolant $g(\tilde{\bsy{\mu}})$ of the error estimator is obtained over the whole parameter domain employing only function evaluations in~(\ref{eq:rbf-interpolant}). Thus the learned surrogate of the error estimator is $g(\tilde{\bsy{\mu}})$. It is not difficult to see that computing the error estimator over the parameter domain is more expensive than using the surrogate $g(\tilde{\bsy{\mu}})$.
\begin{itemize}
\item The cost of computing the surrogate $g(\tilde{\bsy{\mu}})$ for {\it all the} parameter samples $\tilde{\bsy{\mu}}$ in a certain parameter set with cardinality $N_f$ is:
	\begin{itemize}
				\item Solving small, dense ROM : $\mathcal{O}(r^{3}) \times N_{c}$.
				\item Matrix-vector product to evaluate residual : $\mathcal{O}(nr) \times N_{c}$.
				\item Vector-vector inner product to evaluate \cref{eq:err_est_computable} :
			    $\mathcal{O}(n) \times N_{c}$.
				\item Identify the coefficients $c$ by solving \cref{eq:rbfsystem-saddle}: $\mathcal{O}((N_{c}+D)^{3})$.
				\item Evaluate the interpolant through function evaluation \cref{eq:rbf-interpolant}  over a parameter set with cardinality $N_{f}$:
				$\mathcal{O}((N_{c}+D)) \times N_{f}$.
	\end{itemize}
\item The cost of evaluating the error estimator $\Delta(\tilde{\bsy \mu})$ for {\it all the} parameters samples $\tilde{\bsy{\mu}}$ in a certain parameter set with cardinality $N_f$ are,
	\begin{itemize}
		\item Solving small, dense ROM : $\mathcal{O}(r^{3}) \times N_{f}$.
		\item Matrix-vector product to evaluate residual : $\mathcal{O}(nr) \times N_{f}$.
		\item Vector-vector inner product to evaluate \cref{eq:err_est_computable} :
		$\mathcal{O}(n) \times N_{f}$.
	\end{itemize}	
\end{itemize} 
Here, $n$ is the full order dimension of the system; the reduced size $r$ is as small as $N_{c}+D$, i.e. $r \approx N_{c}+D$. For $N_{f} \gg N_{c}$, it is clear that computing the error estimator is more expensive than computing the surrogate.
\subsection{Adaptive Choice of Interpolation Points with Surrogate Error Estimator}
In \cref{alg:2}, we present the proposed adaptive method to choose interpolation points using a surrogate error estimator. We call the algorithm \textsf{IPSUE} - \textsf{Interpolation Points using SUrrogate error Estimator}.  To follow the learning process in Subsection~\ref{subsec:learn} in practice, we do not consider the entire domain $\R^{d}$, but a fine representation of it given by $\Xi_{f} := [\tilde{\bsy{\mu}}_{1}, \ldots, \tilde{\bsy{\mu}}_{N_{f}}]^{\textsf{T}}$, containing $N_{f} \gg N_{c}$ parameters. Therefore,
we consider two training sets: a coarse training set $\Xi_{c}$ and a fine training set $\Xi_{f}$. In Step 4 of \cref{alg:2}, we perform Steps 4 - 6 from \cref{alg:1}. In Step 5, the error estimator is evaluated only over the coarse training set, an important distinction from \cref{alg:1}. In Step 6, the argument of the maximum is chosen as the next interpolation point for $\bV$. Step 7 selects the parameter that maximizes the second summand of the error estimator and uses it as the interpolation point ($\tilde{\bsy \mu}_{\alpha}^{i}$) for enriching $\bV_{e}$ in the next iteration.  As noted in \cite{morFenB19}, it is important that the interpolation points $\tilde{\bsy{\mu}}^{i}$ and $\tilde{\bsy{\mu}}_{\alpha}^{i}$ are distinct, in order to ensure that $\bV_{\text{du}} \neq \bV_{e}$. 
Then, in Step 8, using $\Delta(\tilde{\bsy{\mu}}) \forall\, \tilde{\bsy{\mu}} \in \Xi_{c}$ we learn the error estimator over the parameter domain (represented by $\Xi_{f}$) by determining $g(\tilde{\bsy{\mu}}), \forall\, \tilde{\bsy{\mu}} \in \Xi_{f}$. In Step 9, $n_{a}^{(1)}$ new parameters are identified from $\Xi_{f}$ such that they have the largest errors measured by $g(\tilde{\bsy{\mu}})$. The coarse training set is then updated with the newly identified points.
%
% Algorithm 1 -- 
\renewcommand{\algorithmicrequire}{\textbf{Input:}}
\renewcommand{\algorithmicensure}{\textbf{Output:}}
\begin{algorithm}[t!]
	\caption{\textsf{Interpolation Points using SUrrogate error Estimator (IPSUE)} algorithm}\label{alg:IPSUE}
	\begin{algorithmic}[1]
		\label{alg:2}
		\REQUIRE System matrices $\bA(\bsy{\mu}), \bB(\bsy{\mu}), \bC(\bsy{\mu})\,$, coarse training set $\Xi_{c}$ of cardinality $N_{c}\,$, fine training set $\Xi_{f}$ of cardinality $N_{f}\,$ covering the interesting parameter ranges, tolerance $\epsilon_{tol}$.

		\ENSURE Projection matrix $\bV$.
		
		\STATE Initialize $\bV = [\,]$, $\bV_{\text{du}} = [\,]$,  $\bV_{e} = [\,]$, $\epsilon = 1 + \epsilon_{tol}$, fix $\eta$, the number of moments to be matched. Set $i = 1$.
		
		\STATE Initial interpolation point $\tilde{\bsy{\mu}}^{1}$: a random sample from $\Xi_{c}$ selected using \texttt{rand}(), $\tilde{\bsy{\mu}}_{\alpha}^{1}$: another random sample from $\Xi_{c}$ selected using \texttt{rand}(). Here, \texttt{rand}() is the intrinsic MATLAB\textsuperscript \textregistered function.
		
		\WHILE{$\epsilon$ $>$ $\epsilon_{tol}$} 
		
		\STATE Perform Steps 4 - 6 from \cref{alg:1}.	
		
		\STATE Use \cref{eq:err_est_computable} and obtain $\Delta(\tilde{\bsy{\mu}})\,\, \forall, \tilde{\bsy{\mu}} \in \Xi_{c}$.
		
		\STATE $\tilde{\bsy{\mu}}^{i+1} = \arg \max \limits_{\tilde{\bsy{\mu}} \in \Xi_{c}}  \Delta(\tilde{\bsy{\mu}})$.

%		\STATE $\tilde{\bsy{\mu}}_{\alpha}^{i} = \arg \max \limits_{\tilde{\bsy{\mu}} \in \Xi_{f}} s_{2}(\tilde{\bsy{\mu}})$ (optional).

		\STATE $\tilde{\bsy{\mu}}_{\alpha}^{i+1} = \arg \max \limits_{\tilde{\bsy{\mu}} \in \Xi_{c}} \big|\tilde{\mathbf{e}}_{\text{du}}^{\textsf{T}}(\tilde{\bsy{\mu}}) \br_{\text{pr}}(\tilde{\bsy{\mu}})\big|$.		
		
		\STATE Form the RBF interpolant $g(\tilde{\bsy{\mu}})$ of the error estimator $\Delta(\tilde{\bsy{\mu}})$ over $\Xi_{f}$.
		
%		\STATE Form the RBF interpolant $s_{2}(\tilde{\bsy{\mu}})$ of the second summand of $\Delta(\tilde{\bsy{\mu}})$ over $\Xi_{f}$.
		
		\STATE Select $n_{a}^{(1)}$. Identify  $\big( \tilde{\bsy{\mu}}_{1}^{(1)}, \ldots, \tilde{\bsy{\mu}}_{n_{a}}^{(1)} \big)$ from $\Xi_{f}$ with the largest errors for $g(\tilde{\bsy{\mu}})$. Usually, $n_{a}^{(1)} = 1$.
		
		\STATE Update the coarse training set with the newly identified parameters,\\ $\Xi_{c} := \big[ \Xi_{c} \cup \,\,  \big( \tilde{\bsy{\mu}}_{1}^{(1)}, \ldots, \tilde{\bsy{\mu}}_{n_{a}}^{(1)} \big) \big]$.
			
		\STATE $i = i + 1$. 
		
		\STATE $\epsilon = \Delta(\tilde{\bsy{\mu}}^{i})$.	
		\ENDWHILE
	\end{algorithmic}
\end{algorithm}
%%%%%%%%%%%%%%%%%%%%%%%%%%%%%%%%%%%%%%%%%%%%%%%%%%%%%%%%%%%%%%%%%%%%%%%%%%%%%
\section{Numerical examples}
	\label{sec:4}
	In this section, we provide numerical results to show the efficiency of the proposed \textsf{IPSUE} algorithm. The first example is from circuit simulation used in \cite{morFenB19}. It is characterized by its large parameter range. The second is a benchmark example of a microthruster device, from the \texttt{MORwiki} collection \cite{morWiki}. This model has $4$ parameters. The final example is a finite element model of a waveguide filter, from \cite{Rub18}. All numerical tests were performed in MATLAB\textsuperscript \textregistered 2015a, on a laptop with Intel\textsuperscript \textregistered Core\textsuperscript \texttrademark i5-7200U @ 2.5 GHZ, with 8 GB of RAM. In the numerical results, $N_{\mu}$ refers to the cardinality of the fixed training set $\Xi$, used in \cref{alg:1}, $N_{c}, N_{f}$ are, respectively, the cardinality of the coarse and fine training sets used in \cref{alg:2} and finally $N_{t}$ denotes the cardinality of the parameter test set $\Xi_{t}$ used for validating the accuracy of the final ROMs constructed by \cref{alg:1,alg:2}. Also, we use the same test sets for comparing the performances of \cref{alg:1,alg:2}. 
	\subsection{RLC Interconnect Circuit}
\begin{figure}[t!]
			\centering
			\includegraphics[scale = 0.7]{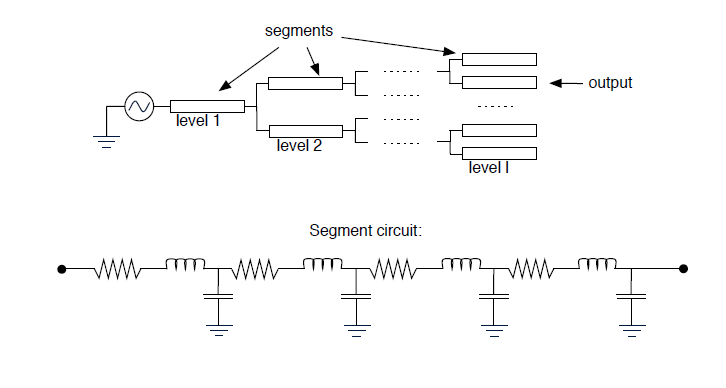}
			\caption{RLC interconnect circuit.}
			\label{fig:RLC}
		\end{figure}

This example models the large-scale interconnects in integrated circuit (IC) design. It is represented in~\cref{fig:RLC}.
			\begin{table}[t!]
				\begin{center}
				\begin{tabular}{|c|c|}
					\hline
					Setting            &   Value \\ \hline \hline
					$n$ 			   &   $6134$ \\ \hline
					$\epsilon_{tol}$   &   $10^{-3}$ \\ \hline
					$N_{\mu}$          &   $90$ parameters \\ \hline
					$N_{c}$            &   $\{ 21, 27 \}$ parameters \\ \hline
					$N_{f}$            &   $200$ parameters \\ \hline
					$N_{t}$            &   $900$ paramters \\ \hline
					$\eta$             &   $3$ \\ \hline
				\end{tabular}
				\end{center}
				\caption{Simulation settings for the RLC interconnect circuit.}
				\label{tab:RLC}
			\end{table}		
		The discretized model has dimension $n = 6134$. It is a non-parametric system in time domain, but in the frequency domain, the frequency $f$ is considered as the parameter and the interpolation points are selected from a wide frequency range: $f \in [0 \,, 3]$GHz. \cref{tab:RLC} gives the simulation settings used for implementing \cref{alg:1,alg:2} to generate the reduced order models for this example.
			\paragraph{\textit{Test 1: \cref{alg:1} applied to RLC model}}
To enable comparison, we use the same training set $\Xi$ used in \cite{morFenB19}. It consists of 90 samples covering the range of interest. The sampled frequencies are given by $f_{i} = 3 \times 10^{i/10}, \,\, s_{i} = 2 \pi \jmath f_{i}$ with $i = 1, 2, \ldots, 90$. \cref{alg:1} converges to the set tolerance in just $3$ iterations. The obtained ROM is of dimension $r = 20$. On average, it takes $3.3$ seconds for the algorithm to converge. For the sake of robustness, we test the ROM on a different set of test parameters $\Xi_{t}$ with $N_{t} = 900$ parameters. \cref{fig:RLC_testset_alg1} shows the error of $\hat{\bH}(s)$ for the parameters in $\Xi_{t}$.
				\newlength\fheight
				\newlength\fwidth 
				\begin{figure}[t]
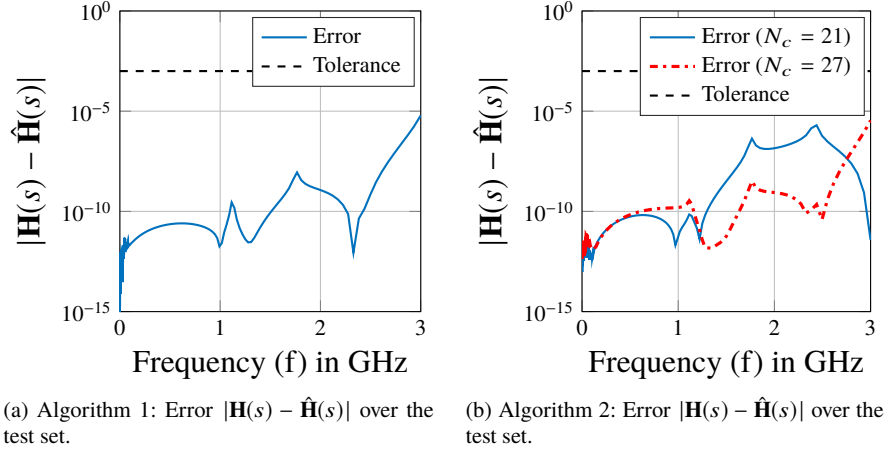

					\centering		
					\setlength\fheight{4cm}
					\setlength\fwidth{4cm}
					\subfloat[\cref{alg:1}: Error $| \bH(s) - \hat{\bH}(s) |$ over the test set.]{\label{fig:RLC_testset_alg1}\input{figures/RLC_testset_trueerror.tex}} \hfill
					\subfloat[\cref{alg:2}: Error $| \bH(s) - \hat{\bH}(s) |$ over the test set.]{\label{fig:RLC_testset_alg2}\input{figures/RLC_testset_trueerror_Test2.tex}}							
					\caption{Results for the RLC  model.}
					\label{fig:RLC_testset_error}
				\end{figure}
			\paragraph{\textit{Test 2: \cref{alg:2} applied to RLC model}}
				Next, we test \cref{alg:2} on the RLC interconnect model. For this, we consider two different coarse training sets $\Xi_{c}$ of cardinality $21, 27$ sampled as, $\Xi_{c}^{j} = 3 \times 10^{j/10},\,\, j = 1, 2, \ldots, 21$ and $j = 1, 2, \ldots, 27$. We consider different samplings for the fine training set in order to numerically illustrate that the proposed algorithm is independent of the kind of sampling used. The fine training set $\Xi_{f}$ consists of $200$ logarithmically distributed parameters in the first case and in the second case contains 200 parameters distributed as $\Xi_{f}^{j} = 3 \times 10^{j/10},\,\, j = 1, 2, \ldots, 200$. For the RBF interpolation, we use thin-plate splines as the kernel function. It is given by $\Phi(\tilde{\bsy{\mu}}_{1} , \tilde{\bsy{\mu}}_{2}) := (\| \tilde{\bsy{\mu}}_{1} - \tilde{\bsy{\mu}}_{2} \|_{2})^{2} \log_{e}(\| \tilde{\bsy{\mu}}_{1} - \tilde{\bsy{\mu}}_{2} \|_{2}) $.
				\cref{alg:2} converges to the specified tolerance in just $3$ iterations for both choices of $\Xi_{c}$, with $n_{a}^{(1)} = 1$. In the first case, the obtained ROM is of dimension $r = 21$ and takes $1.6$ seconds to converge in average. The second case results in a ROM of dimension $21$ and takes $1.7$ seconds on average to converge to the defined tolerance. \cref{fig:RLC_testset_alg2} shows the error of $\hat{\bH}(s)$ at parameters in the test set $\Xi_{t}$ produced by the ROM obtained using \cref{alg:2}, with two different coarse training sets. Clearly, \cref{alg:2} takes less time than \cref{alg:1}, while still producing a ROM that is uniformly below the tolerance, on an independent test set.				
				
\subsection{Thermal Model}
The second example is the model of the heat transfer inside a microthruster unit \cite{morWiki}. It is obtained after spatial discretization using the finite element method and has dimension $n = 4257$. The governing equation is given as,
\begin{table}[t!]
	\begin{center}
		\begin{tabular}{|c|c|}
			\hline
			Setting            &   Value \\ \hline \hline
			$n$ 			   &   $4257$ \\ \hline
			$\epsilon_{tol}$   &   $10^{-4}$ \\ \hline
			$N_{\mu}$          &   $625$ parameters, log-sampled \\ \hline
			$N_{c}$            &   $256$ parameters, log-sampled \\ \hline
			$N_{f}$            &   $2401$ parameters, log-sampled \\ \hline
			$N_{t}$            &   $1000$ parameters, randomly sampled\\ \hline
			$\eta$             &   $1$ \\ \hline
		\end{tabular}
	\end{center}
	\caption{Simulation settings for the thermal model.}
	\label{tab:thermal}
\end{table}
	\begin{equation*}
		\begin{aligned}
			\bE \dot{\bx}(t) &= (\bA_{0} - \sum_{i=1}^{3} h_{i} \bA_{i}) \bx(t) + \bB \bu(t),\\
			\by(t) &= \bC \bx. 
		\end{aligned}		
	\end{equation*}
Here, $\bE, \bA_{0}$ are symmetric sparse matrices representing the heat capacity and heat conductivity, respectively. $\bA_{i},\, i \in \{1,2,3\}$ are diagonal matrices governing the boundary condition. The parameters $h_{1}, h_{2}, h_{3} \in [1, \, 10^{4}]$ represent, respectively, the film coefficients of the top, bottom and side of the microthruster unit. We transform the above system to the frequency domain and apply \cref{alg:1} and \cref{alg:2}. In the frequency domain, the system has four parameters $\tilde{\bsy{\mu}} := (s, h_{1}, h_{2}, h_{3})$ with $s = \jmath 2 \pi f$. The frequency range of interest is $f \in [10^{-2},\,10^2]$ Hz. The tolerance for the ROM is set as $10^{-4}$.
				\begin{figure}[t]
					\centering		
					\setlength\fheight{4cm}
					\setlength\fwidth{4cm}
					\subfloat[\cref{alg:1}: Error $| \bH(\tilde{\bsy \mu}) - \hat{\bH}((\tilde{\bsy \mu}) |$ over the test set.]{\label{fig:thermal_testset_alg1}\input{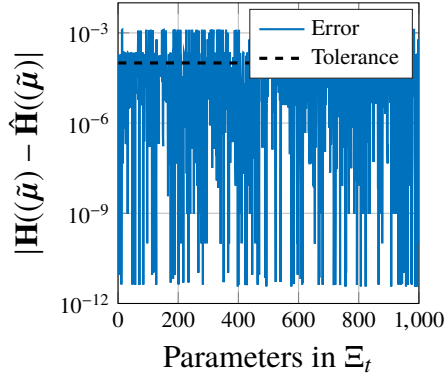}} \hfill
					\subfloat[\cref{alg:2}: Error $| \bH((\tilde{\bsy \mu}) - \hat{\bH}((\tilde{\bsy \mu}) |$ over the test set.]{\label{fig:thermal_testset_alg2}\input{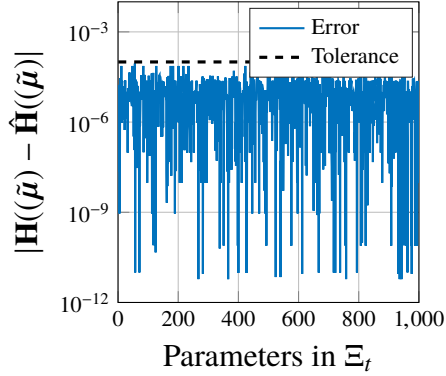}}							
					\caption{Results for the thermal  model.}
					\label{fig:thermal_testset_error}
				\end{figure}
			\paragraph{\textit{Test 3: \cref{alg:1} applied to the thermal model}}
Owing to the wide range of parameters, we consider a large fixed training set ($\Xi$). To construct it, we consider 5 logarithmically-spaced samples for each of the 4 parameters and form a grid consisting of $5^4$ samples. For the test set $\Xi_{t}$, we form a grid of $8^{4}$ logarithmically spaced parameters and randomly select $1000$ parameters from it. The greedy algorithm takes 10 iterations to converge and results in a ROM of size $r = 86$. On an average over 5 runs, the greedy algorithm takes 254 seconds to converge.
In \cref{fig:thermal_testset_alg1}, we see the performance of the resulting ROM over $\Xi_{t}$. For several parameters, the ROM fails to meet the desired tolerance. This indicates that the training set was not fine enough to capture all the variations in the solutions over the parameter domain.
			\paragraph{\textit{Test 4: \cref{alg:2} applied to the thermal model}}
We now consider \cref{alg:2} applied to the thermal model. The coarse training set $\Xi_{c}$ has $4^{4}$ parameters, with logarithmic sampling. The fine training set $\Xi_{f}$ has $7^{4}$ parameters. For the RBF interpolation, we make use of thin-plate splines as the kernel function. Further, we set $n_{a}^{(1)} = 1$ in Step 9 of \cref{alg:2} so that the coarse training set is updated with one new parameter per iteration. The same test set as in Test 3 is used. The resulting ROM has order $r = 85$ and its error over the test set is below the tolerance, as shown in \cref{fig:thermal_testset_alg2}. The algorithm took 162 seconds to converge. The runtime was measured as an average over 5 independent runs of the algorithm. 
Compared with \cref{alg:1}, \cref{alg:2} is able to meet the required tolerance with a much smaller training set and also in shorter time.
	\subsection{Dual-Mode Circular Waveguide Filter}
\begin{figure}[t!]
	\centering
	\includegraphics[scale = 0.35]{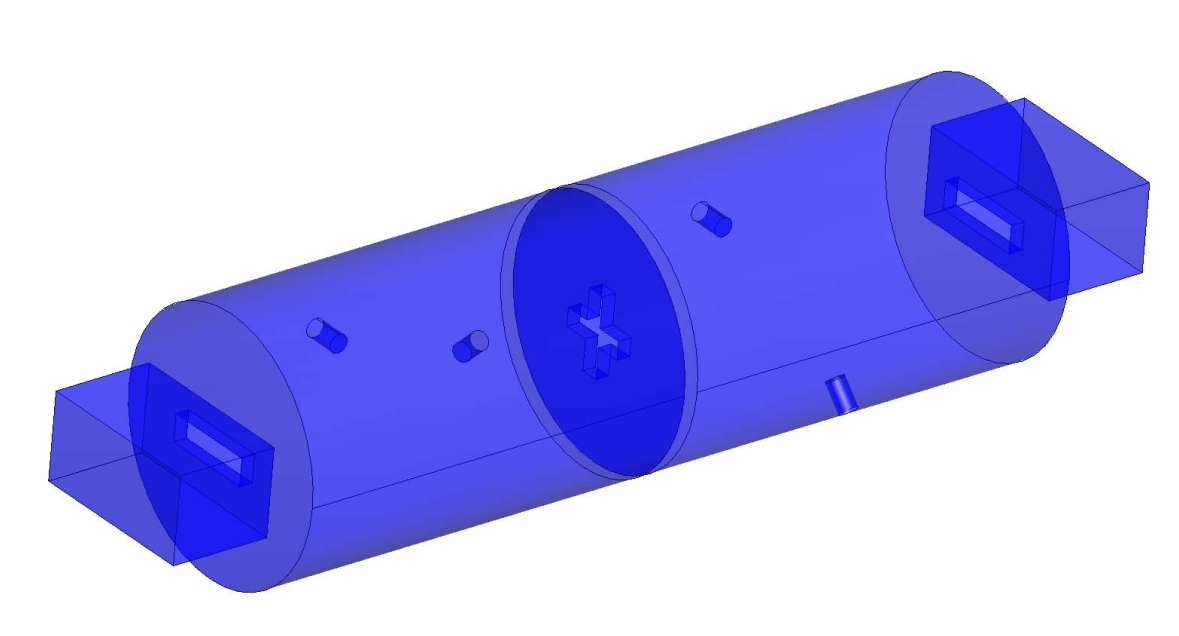}
	\caption{Dual-mode waveguide filter model from \cite{Rub18}.}
	\label{fig:filter}
\end{figure}

The next example is a MIMO system based on the model of a dual-mode circular waveguide filter from \cite{Rub18}, see \cref{fig:filter}. It is a type of narrow bandpass filter widely used in satellite communication due to its power handling capabilities.
Its operation is governed by the time-harmonic Maxwell's Equations. After discretization in space, the governing equations of the filter can be represented in the form of \cref{eq:primal_state}. The system consists of just the frequency parameter $s := \jmath 2 \pi f$, where $f \in [11.5\,,12]$ GHz is the operating frequency band of the filter. The affine form of the system matrix is,  $\mathscr{A}(s) := \mathcal{S} + s^{2} \mathcal{T}$ and $\bB(s) := s \mathcal{Q}$. We have $\mathcal{S}, \mathcal{T} \in \R^{n \times n}$ and $\mathcal{Q} \in \R^{n \times 2}$, with $n = 36426$. The system has two inputs and two outputs. \cref{tab:filter} summarizes the simulation settings.
The quantity of interest are the scattering parameters, obtained via post-processing \cite{RubRM09} from the system output $\by(s) := \mathcal{Q}^{\textsf{T}} \bX(s)$. It is easy to see that $\by(s)$ has the same expression as $\bH(\tilde{\bsy \mu})$ in \cref{eq:LTI_s} for $\tilde{\bsy \mu}=s$. The error estimator $\Delta(\tilde{\bsy \mu})$ in \cref{sec:greedy} can be directly applied to estimate the error of $\hat{\by}(s)$ computed by the ROM. See \cite{morFenB19} for detailed analysis.
Since the system is MIMO, the scattering parameters at a given $s$ are in the form of a complex-valued matrix given by
	\begin{equation*}
		\bS  = \begin{bmatrix}
		S_{11} & S_{12} \\
		S_{21} & S_{22}
		\end{bmatrix}.
	\end{equation*}
Scattering parameters are important in characterizing the performance of filters \cite{pozar1998microwave}.
				\begin{figure}[t]
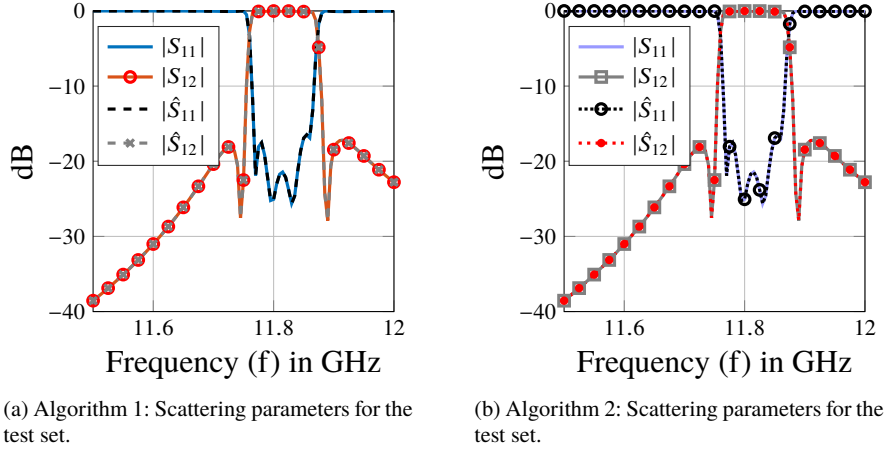

	\centering		
	\setlength\fheight{4cm}
	\setlength\fwidth{4cm}
	\subfloat[\cref{alg:1}: Scattering parameters for the test set.]{\label{fig:filter_testset_sparam_alg1}\input{figures/filter_fixedset_sparam.tex}} \hfill
	\subfloat[\cref{alg:2}: Scattering parameters for the test set.]{\label{fig:filter_testset_sparam_alg2}\input{figures/filter_adaptsetset_sparam.tex}}						
	\caption{Scattering parameters for the dual-mode filter.}
	\label{fig:filter_sparam}
\end{figure}
\begin{figure}[t]
	\centering		
	\setlength\fheight{4cm}
	\setlength\fwidth{4cm}
	\subfloat[\cref{alg:1}: Error of the scattering parameters computed from the ROM over the test set.]{\label{fig:filter_fixedset_trueerror_alg1}\input{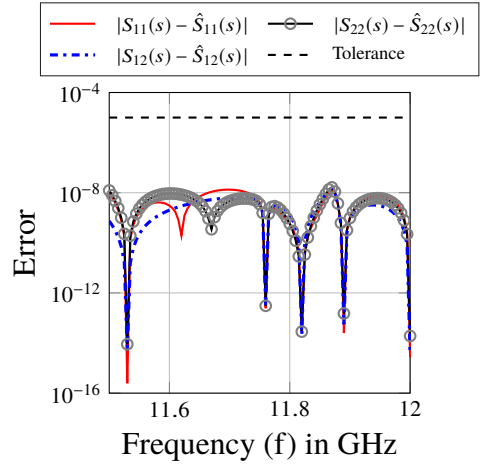}} \hfill
	\subfloat[\cref{alg:2}: Error of the scattering parameters computed from the ROM over the test set.]{\label{fig:filter_adaptset_trueerror_alg2}\input{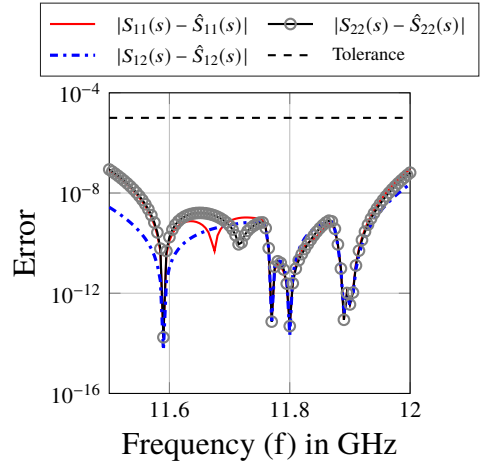}}						
	\caption{Results for the dual-mode filter.}
	\label{fig:filter_trueerror}
\end{figure}
\begin{table}[t!]
	\begin{center}
		\begin{tabular}{|c|c|}
			\hline
			Setting            &   Value \\ \hline \hline
			$n$ 			   &   $36426$ \\ \hline
			$\epsilon_{tol}$   &   $10^{-5}$ \\ \hline
			$N_{\mu}$          &   $51$ parameters, uniformly spaced \\ \hline
			$N_{c}$            &   $17$ parameters, uniformly spaced \\ \hline
			$N_{f}$            &   $500$ parameters, sampled randomly \\ \hline
			$N_{t}$            &   $101$ parameters, uniformly spaced\\ \hline
			$\eta$             &   $1$ \\ \hline
		\end{tabular}
	\end{center}
	\caption{Simulation settings for the dual-mode filter.}
	\label{tab:filter}
\end{table}
			\paragraph{\textit{Test 5: \cref{alg:1} applied to the dual-mode filter}}
Applying \cref{alg:1} with the fixed training set $\Xi$ to the model results in a ROM of size $r = 10$ with the greedy algorithm taking $5$ iterations to converge. Since this example is a MIMO system, we make use of \cref{eq:error_estm_MIMO}. The average runtime over 5 independent runs of \cref{alg:1} was found to be 46 seconds.
\cref{fig:filter_testset_sparam_alg1} plots the scattering parameters computed from FOM simulations and those obtained from the ROM at the parameters in the test set. We plot the absolute values of full order scattering parameters $S_{11}, S_{12}$ and the corresponding reduced ones $\hat{S}_{11}, \hat{S}_{12}$ on a Decibel scale. 		
In \cref{fig:filter_fixedset_trueerror_alg1} we plot the error of the scattering parameters $\hat{S}_{11}, \hat{S}_{12}, \hat{S}_{22}$ computed from the ROM, over the test set $\Xi_{t}$. Note that since $S_{12} = S_{21}$, we only show the error $|S_{12} - \hat{S}_{12}|$.
			\paragraph{\textit{Test 6: \cref{alg:2} applied to the Dual-mode filter}}
In Step 8 of \cref{alg:2}, we construct an RBF surrogate for each of $\Delta_{ij},\, i, j \in \{1, 2\}$ in \cref{eq:error_estm_MIMO} for this MIMO system.  $n_{a}^{(1)}$ is set to be 1 and inverse multiquadric is used as the kernel function. It is given by $\Phi(\tilde{\bsy{\mu}}_{1}, \tilde{\bsy{\mu}}_{2}) := \frac{1}{1 + \big(\gamma \| \tilde{\bsy{\mu}}_{1} - \tilde{\bsy{\mu}}_{2}  \|_{2} \big)^{2} } $. $\gamma$ is a user-defined parameter and we set $\gamma = 16$ in our experiments.
We pick the maximum among the four surrogates and add the corresponding parameter to the coarse training set, i.e. in Step 9 of \cref{alg:2}, we replace $g(\tilde{\bsy \mu})$ with $\max \limits_{i,j \in \{1,2\}} g_{i,j}(\tilde{\bsy \mu})$.
\cref{alg:2} results in a ROM that is of the same size as the ROM from Test 5 ($r = 10$). However, on average, \cref{alg:2} only needs $24$ seconds to converge, almost half that of the time required in Test 5. In \cref{fig:filter_adaptset_trueerror_alg2} we plot the  errors of the scattering parameters computed from the ROM over the test set $\Xi_{t}$. \cref{fig:filter_testset_sparam_alg2} plots the scattering parameters from the FOM simulations and those computed by the ROM. Both algorithms result in ROMs meeting the specified tolerance, but \cref{alg:2} requires much shorter time to generate the ROM.
%		
%%%%%%%%%%%%%%%%%%%%%%%%%%%%%%%%%%%%%%%%%%%%%%%%%%%%%%%%%%%%%%%%%%%%%%%%%%%
\section{Conclusion}
	\label{sec:5}	
	In this work, we have proposed \textsf{IPSUE}, an adaptive algorithm for updating the training set and choosing the interpolation points for frequency domain MOR methods. Our target applications are cases where the problem parameters vary over a wide range of values, or the parameter space dimension is larger than two. In either of these cases, many interpolatory MOR algorithms may take a long time to generate the ROM. Moreover, a naive, heuristic sampling of the parameter training set may result in a ROM that is not robust. \textsf{IPSUE} offers a viable means to generate reliable ROMs that satisfy the user-defined tolerance and at the same time, without being offline expensive. The illustrated numerical examples show that it is a promising approach. As future work, we plan to apply the algorithm to more complex models.
\begin{acknowledgement}
The first author is affiliated to the International Max Planck Research School for Advanced Methods in Process and Systems Engineering (IMPRS-ProEng).	
\end{acknowledgement}	
 \bibliographystyle{spmpsci}
 \bibliography{myref}
\end{document}